\documentclass{IEEEtran}
\usepackage[utf8]{inputenc}

\usepackage{graphicx}
\usepackage{amsmath}
\usepackage{amssymb}
\usepackage{amsfonts}
\usepackage{mathrsfs}
\usepackage{dsfont}
\makeatletter
\def\amsbb{\use@mathgroup \M@U \symAMSb}
\makeatother
\usepackage{color,cite}
\usepackage{epstopdf}

\usepackage[usenames,dvipsnames]{xcolor}
\newtheorem{theorem}{Theorem}
\newtheorem{remark}{Remark}
\newtheorem{definition}{Definition}
\newtheorem{lemma}{Lemma}
\newtheorem{proposition}{Proposition}

\newtheorem{assumption}{Assumption}

\newcommand{\vect}[1]{\boldsymbol{#1}}

\usepackage{enumerate}
\usepackage{balance}

\usepackage[prependcaption,colorinlistoftodos]{todonotes}
\newcommand{\ak}[1]{{\color{black}{#1}}}
\newcommand{\icl}[1]{{\color{black}{#1}}}

\newcommand{\closure}[2][3]{%
  {}\mkern#1mu\overline{\mkern-#1mu#2}}

\title{\Large \bf{Secondary frequency control with on-off load side participation in power networks}}

\author{Andreas Kasis\thanks{\icl{This work was supported by ERC starting grant 679774.}}
\thanks{Andreas Kasis and Ioannis Lestas are with the Department of Engineering, University of Cambridge, Trumpington Street, Cambridge, CB2 1PZ, United Kingdom; e-mails: ak647@cam.ac.uk, icl20@cam.ac.uk}, Nima Monshizadeh\thanks{Nima Monshizadeh is with the Engineering and Technology Institute, University of Groningen, Nijenborgh 4, 9747AG, Groningen, The Netherlands.
email: n.monshizadeh@rug.nl} and Ioannis Lestas
\thanks{A preliminary version of this work {has appeared} in \cite{kasis2017switch_CDC}. This manuscript includes the proofs of the main results, {generalizations to wider classes of systems}, {as well as} additional results and discussion {that demonstrate} the applicability of the proposed analysis.}}

\begin{document}

\maketitle

\begin{abstract}
We study the problem of {decentralized} secondary frequency regulation in power networks where ancillary services are provided via on-off load-side participation.
We initially consider on-off loads that switch when prescribed frequency thresholds are exceeded, together with a large class of passive continuous dynamics for generation and demand.
The considered on-off loads are able to assist existing secondary frequency control mechanisms and return to their nominal operation when the power system is restored to its normal operation, a highly desirable feature which minimizes users disruption. We show that system stability is not compromised despite the switching nature of the loads. However, such control policies are prone to {chattering}, which limits the practicality of these schemes. As a remedy to this problem, we propose a {hysteretic on-off policy where loads switch on and off at different frequency thresholds and show that stability guarantees are retained when the same decentralized passivity conditions for continuous generation and demand hold.
Several relevant examples
are discussed to demonstrate the applicability of the proposed results. Furthermore,}
we verify our analytic results with numerical investigations on the Northeast Power Coordinating Council (NPCC) 140-bus system.
\end{abstract}

\section{Introduction}

\textbf{Motivation and literature survey:} It is anticipated that renewable sources of generation will increase their penetration
in power networks in the near future \cite{lund2006large, ipakchi2009grid}. This is expected to introduce intermittency in the power generated resulting {in additional challenges in the real time operation of power networks that need to be addressed.}

{A main objective in the operation of a power system is
to ensure that generation matches demand in real time. This is achieved by means of primary and secondary  frequency control schemes with the latter also ensuring that the frequency returns {to} its nominal value (50Hz or 60Hz).}
 Secondary frequency control {is} traditionally performed by having the generation side following demand \cite{kundur1994power}. However, a large penetration of renewable sources of generation limits the controllability of generation and at the same time {makes} the system more sensitive to disturbances due to the reduced system inertia \cite{ulbig2014impact}.
Controllable loads are considered by many {a promising} solution to counterbalance intermittent generation, being able to adapt their demand based on frequency deviations,  providing fast response at urgencies.
 Recently, various research studies focused on the inclusion of controllable demand to aid both primary control as in \cite{molina2011decentralized, trip2016optimal, zhaotopcu, kasis2016primary, devane2016primary, kasis2017primary}   and secondary control  as in \cite{P30_Low, trip2016internal,kasis2017stability}.

Further from providing ancillary services at urgencies, it is also desired that controllable loads are non-disruptive, i.e. their  assistance {should have} a negligible effect on users comfort, {see e.g. \cite{liu2016non}. Non-disruptive load-side control schemes ensure that}
loads alter their demand at urgencies but return to their normal operation when the danger for the network has been surpassed.
Moreover, in many occasions, a realistic representation of loads involves only a discrete set of  possible demand values, e.g. on and off states.
Hence,  {incorporating on-off controllable loads that appropriately react to frequency deviations in power networks is of particular interest in load-side participation schemes.}

\textbf{Contribution:} As already mentioned, the incorporation of controllable demand within the power network has been extensively studied in literature \cite{molina2011decentralized, trip2016optimal, zhaotopcu, kasis2016primary, devane2016primary, kasis2017primary, P30_Low,trip2016internal, kasis2017stability}. However, {existing theoretical studies consider {power networks with} controllable loads that vary continuously with deviations in frequency. These are thus unable to {capture} discontinuities, which are relevant when loads can take only a discrete number of states, e.g. on and off, as is often the case in practice.
This study provides analytical results which allow on-off loads to contribute to frequency control within a power network with a general network topology, when the continuous part of the generation/demand dynamics, that are allowed to be heterogeneous and of higher order, satisfy certain decentralized passivity conditions. In particular,} by the prospect of loads providing ancillary service to the power network in a {decentralized} way and the desire that users comfort will be distorted for only short periods of time, we consider controllable on-off loads that switch when some frequency deviation is reached so that they assist the network at urgencies (i.e. when large frequency deviations are experienced) and {otherwise return to their original operation.}
{It will be shown that the inclusion} of such loads does not compromise the stability of the power network, and {results} in enhanced frequency performance. However, it will be observed {that such controllable loads may switch arbitrarily fast} within a finite interval of time, or in other words, {exhibit} {chattering}. To avoid {this {phenomenon},} we propose on-off loads with hysteretic dynamics {where loads switch on and off at different frequency thresholds.
{{As it will be analytically shown in the paper}, unlike (instantaneous) on-off loads, hysteretic loads do not induce {chattering} behavior.}
{Furthermore, we show that
stability guarantees can be provided {when such schemes are used within the power network.}
{It should be noted that hysteretic dynamics increase the complexity of the analysis as they render the underlying dynamical system a hybrid system, thus requiring appropriate analysis tools to be employed \cite{goebel2012hybrid}.
Various examples are provided to demonstrate the applicability of the {results presented}.
Moreover,} we provide a numerical validation of our analytic results on the NPCC 140-bus system, {where it is demonstrated that the incorporation of frequency dependent on-off loads provides improved frequency response by reducing its overshoot.}

\textbf{Paper structure:} The structure of the paper is as follows. Section \ref{Notation} includes some basic notation and preliminaries and in sections \ref{sec:Network_model} and \ref{Section: Passive_continuous} we present the power network model and the considered continuous generation and demand dynamics respectively. In section \ref{Sec:switch} we {study the behavior of {controllable demand that {switches} on/off whenever certain frequency thresholds are met and present} our results concerning network stability. {We also discuss in this section that problematic {chattering} behavior can be observed.}
In section \ref{Sec:Hysteresis}, {we consider controllable loads with hysteretic patterns and} show that {the} stability results extend to this case {while avoiding {chattering}}.}
{In section \ref{Sec: Discussion} we provide relevant examples of continuous generation and demand dynamics, including higher order schemes, that fit within the proposed framework and allow for stability guarantees to be deduced when {discontinuous and hysteretic} loads are also considered.
{Numerical investigations of the results are provided  in section \ref{Simulation}}.
Finally, conclusions are drawn in section \ref{Conclusion}. {The appendix includes the proofs of the results.}}

\section{Notation}\label{Notation}

{Real numbers are denoted by $\mathbb{R}$, and the set of n-dimensional vectors with real entries is denoted by $\mathbb{R}^n$.
{The set of non-negative real numbers is denoted by $\mathbb{R}_{\geq 0}$.}
The set of natural numbers, including zero, is denoted by $\mathbb{N}_0$.} We use $\vect{0}_n$ to denote  the $n \times 1$ vector with all elements equal to $0$.
{A function $f:\mathbb{R}^n \rightarrow \mathbb{R}$
is said to be positive definite
if $f(0) = 0$  and $f(x) > 0$ for every non-zero $x \in \mathbb{R}^n$. It is positive semidefinite if the inequality $>0$ is replaced by
$\ge0$.}
For a discrete set $\Sigma$, let $|\Sigma|$ denote its cardinality.
{For a scalar valued function $V(x), \ V:\mathbb{R}^n\rightarrow \mathbb{R}$, we denote its gradient by $\frac{\partial V}{\partial x}$. We denote by {${\rm sgn}(x)$} the function that takes the} value of $1$ when $x$ is non-negative and $-1$ otherwise. {The Laplace transform of a signal $h(t)$
 is denoted by $\hat{h}(s) = \int_0^\infty e^{-st} h(t) \, dt$.}
{A set {$M\subset{\mathbb{R}^n}$} is said to be invariant {with respect to a system $\dot{x} = f(x)$, where $f:\mathbb{R}^n \rightarrow \mathbb{R}^n$ is locally Lipschitz, if $x(0) \in M$ implies \ak{that the solution $x(t)$ satisfies} $x(t) \in M, \forall t \in \mathbb{R}$.}
Finally, the closure of a set $S$ is denoted by $\closure{S}$. }

\section{Network model}\label{sec:Network_model}

We describe the power network model by a connected graph $(N,E)$ where $N = \{1,2,\dots,|N|\}$ is the set of buses and $E \subseteq N \times N$ the set of transmission lines connecting the buses.
Furthermore, we use $(i,j)$ to denote the link connecting buses $i$ and $j$ and assume that the graph $(N,E)$ is directed with arbitrary {orientation}, so that if $(i,j) \in E$ then $(j,i) \notin E$. For each $j \in N$, we use $i:i\rightarrow j$ and $k:j\rightarrow k$ to denote the sets of buses that are predecessors and successors of bus $j$ respectively.
The following assumptions are made for the network: \newline
1) Bus voltage magnitudes are $|V_j| = 1$ p.u. for each $j \in N$. \newline
2) Lines $(i,j) \in E$ are lossless and characterized by their susceptances $B_{ij} = B_{ji} > 0$. \newline
3) Reactive power flows do not affect bus voltage phase angles and frequencies.

{
\begin{remark}\label{r:assumptions}
{The above assumptions are standard in secondary frequency control studies, e.g. \cite{molina2011decentralized, jiang2012delay, alomoush2010load, li2016connecting} and are valid when medium to high voltages are considered \cite{Bergen_Vittal} {or tight voltage control is {present.}} {In Section \ref{Simulation}, we also} {demonstrate}
our results {in} a realistic simulation {study.}
}
\end{remark}
}

{We use swing equations, see e.g. \cite{Bergen_Vittal}, to describe the frequency dynamics at each bus:}
\begin{subequations} \label{sys1}
\begin{gather}
\dot{\eta}_{ij} = \omega_i - \omega_j, \; (i,j) \in E, \label{sys1a} \\
\hspace{-2.0mm}M_j \dot{\omega}_j\hspace{-1mm} = \hspace{-1mm}- p_j ^L + p_j^M \hspace{-1.0mm} - \hspace{-1.0mm} (d^c_j + d^u_j) - \hspace{-1.5mm} \sum_{k:j\rightarrow k}\hspace{-1.5mm} p_{jk} + \hspace{-1.5mm}\sum_{i:i\rightarrow j}\hspace{-1.5mm} p_{ij},  j\in \hspace{-1.0mm}N,\hspace{-1.0mm} \label{sys1b} \\
\vspace{-0.3cm}
{p_{ij}=B_{ij} \sin\eta _{ij},} \; (i,j) \in E. \label{sys1d}
\end{gather}
\end{subequations}

{In system~\eqref{sys1}, the time-dependent variables  $p^M_j$, $\omega_j$ represent, respectively,  the mechanical power injection, and the deviation from the nominal\footnote{The nominal value of frequency is equal to 50Hz (or 60Hz)} frequency at bus $j$.
{The time-dependent variable $d^c_j$ represents the deviation of the  frequency-dependent controllable load at bus $j$ from its nominal\footnote{{In particular,  $d^c_j(t)=\tilde{d}^c_j(t)-d_j^{nom}$ where $\tilde{d}^c_j(t)$ is the controllable load and $d_j^{nom}$ is a constant nominal value chosen by the users. To simplify the presentation $d_j^{nom}$ is included within the constant $p^L_j$.}}
value.}
The quantity $d^u_j$ is also a time-dependent variable that represents the uncontrollable frequency-dependent load and generation damping present at bus $j$.
Furthermore,
the quantities $\eta_{ij}$ and $p_{ij}$ are time-dependent variables that represent, respectively, the power angle difference,
and the power transmitted from bus $i$ to bus~$j$.
The constant $M_j > 0$ denotes the generator inertia,  {and the constant $p_j^L$ {includes} the frequency-independent load {and the nominal value of the controllable load} at bus $j$.}
We study the response of system~\eqref{sys1} at a step change in
{the {frequency-independent} uncontrollable demand $p_j^L$ at each bus.}

\section{Passive continuous dynamics}\label{Section: Passive_continuous}

{Before we study the effect {of on-off} loads in the power network, we consider, in this section,} a general class of continuous passive nonlinear dynamics for generation and demand.   Existing studies,  e.g. \cite{kasis2016primary},\cite{kasis2016stability}, demonstrated stability {of the power network when such dynamics are present.}
 In the following sections, we demonstrate that the generation and demand dynamics presented {in} this section allow the incorporation of on-off loads without compromising stability, {despite their non-smooth and hybrid nature}. The {practical relevance and generality of the class of dynamics described in this section} will be demonstrated with various examples in Section \ref{Sec: Discussion}.

\vspace{-0mm}
\subsection{Dynamics for generation and demand}\label{Sec: Passive_dynamics}
 For convenience in the analysis, we define a net power supply variable $s_j$ that represents the aggregation of {generation and uncontrollable {frequency-dependent} demand at bus $j$,} given by
\begin{equation}\label{s_variable}
s_j = p^M_j - d^u_j, j \in N.
\end{equation}
{To incorporate general classes of dynamics for the net supply variables, we assume that $s_j$ is the output of a nonlinear system, namely}
\begin{align}\label{dynsys}
& \dot{x}^s_j = f_j(x^s_j,-\omega_j), \nonumber \\
& s_j = g_j(x^s_j, -\omega_j),
\end{align}
{for each $j\in N$}.
{Here, $x^s_j \in \mathbb{R}^{n_j}$, {is the state} of the system, and  $f_j: \mathbb{R}^{n_j}\times \mathbb{R} \rightarrow
\mathbb{R}^{n_j}$,   and $g_j: \mathbb{R}^{n_j} \times \mathbb{R} \rightarrow \mathbb{R}$ are globally Lipschitz for each $j \in N$.}

{
\subsection{Equilibrium analysis}
{{When $d^{c}(t) = \vect{0}_{|N|}, t \in \mathbb{R}$} system \eqref{sys1}--\eqref{dynsys} represents a classical power system model with {constant demand.}
We assume below that there exists an equilibrium to this system.}

{
\begin{definition} \label{eqbrdef_0}
An equilibrium point of the system~\eqref{sys1}--\eqref{dynsys} with  $d^{c}(t) = \vect{0}_{|N|}, t \in \mathbb{R}$ is a value of the states $(\eta, \omega, x^{s})$ for which the time derivatives
of the system are equal to zero.
\end{definition}
}

{
\begin{assumption}\label{assum_eqbl_gen}
There exists an equilibrium point {$(\eta^*, \omega^*, x^{s,*})$  to {the system} \eqref{sys1}--\eqref{dynsys}, with $d^{c}(t) = \vect{0}_{|N|}, t \in \mathbb{R}$.}
\end{assumption}
}}

{It should be noted that the problem of existence of an equilibrium point to \eqref{sys1}--\eqref{dynsys} with $d^{c}(t) = \vect{0}_{|N|}$ is a problem that has been studied in the literature (e.g. \cite{dorfler2013synchronization}) and sufficient conditions involving large enough line susceptances have been derived.}

{In addition,} we impose a constraint on the differences of the phase angles at the equilibrium. This assumption, stated below, is ubiquitous in power network literature, and is treated as a security constraint.
\begin{assumption} \label{assum1}
$| \eta^*_{ij} | < \tfrac{\pi}{2}$ for all $(i,j) \in E$.
\end{assumption}

{
\subsection{Conditions on generation and demand dynamics}

In this section we present {conditions for the system in \eqref{dynsys}. In particular,}} we make the following assumptions:
{\begin{assumption}\label{assum_0}
{{For each $j\in N$} the system  {$\dot{x}^s_j = f_j(x^s_j, 0)$}} has no compact invariant set other than equilibrium points.
\end{assumption}}
{Assumption \ref{assum_0}
 {implies that when the} frequency deviation in \eqref{dynsys} is zero} then any {compact invariant set {of \eqref{dynsys}}
consists of equilibrium points} only.}
For a
{stable}
linear system, the latter means that the state matrix has {no eigenvalues on the imaginary axis except from a possible single eigenvalue at the origin.}}

Furthermore, the following assumption ensures that at equilibrium, the frequency will be at its nominal value.

\begin{assumption}\label{assum_nominal_freq}
There exists at least one bus $j$ such that {when the {vector field in \eqref{dynsys} satisfies}
{$f_j(\bar{x}^s_j,-\bar{\omega}_j) = 0$}{, for some constants $\bar \omega_j$ and $\bar x_j^s$,} {then} $\bar{\omega}_j = 0$.}
\end{assumption}

\begin{remark}\label{rem_assum2}
Assumption \ref{assum_nominal_freq} requires that for at least one bus, the dynamics are such that an equilibrium can be reached only if the frequency is at its nominal value.
Such requirement is satisfied by control policies where integral action is used and also by more general schemes, as discussed in Section \ref{Sec: Discussion}.
{Moreover, note that Assumption \ref{assum_nominal_freq} needs not to be satisfied by all buses.} This reflects the fact that {secondary} control is performed by a small number of buses, where the rest may only provide support at faster timescales {without incorporating integral action}.
\end{remark}

The final condition on the dynamics in \eqref{dynsys}  is related to a notion of passivity that is discussed next.

\begin{definition}\label{Passivity_Definition}
The system~\eqref{dynsys} is said to be locally input strictly passive {around} the constant {input} $-\bar{\omega}_j$ and the
{point}
~$\bar{x}^s_j$ if there exist open neighborhoods $\Omega_j$ of $-\bar{\omega}_j$ and $X^s_j$ of $\bar{x}^s_j$ and a continuously differentiable positive {semidefinite} function $V^S_j(x^s_j)$ (called the storage function),  such that for all $-\omega_j \in \Omega_j$ and all $x^s_j \in X^s_j$,
\begin{equation}\label{In_passivity_def}
\dot{V}^S_j(x^s_j){ \leq (-\omega_j +\bar{\omega}_j) (s_j - \bar{s}_j) - \phi_j(-\omega_j +\bar{\omega}_j),}
\end{equation}
where $\phi_j$ is a positive definite function and $\bar{s}_j = g_j(\bar{x}^s_j, -\bar{\omega}_j)$.
\end{definition}

We now {consider an equilibrium point {of \eqref{dynsys}} with frequency\footnote{{As discussed later within the paper, it easily follows form Assumption \ref{assum_nominal_freq}
that the equilibrium frequency always satisfies $\omega^*_j=0,$ {$\forall j \in N$}.}}
$\omega^*_j=0$ and equilibrium states $x^{s,*}_j, j\in N$ and
suppose that the power supply variables at each bus satisfy the input strict passivity condition {around} this point.}

\begin{assumption} \label{assum_passivity}
{For each $j\in N$, the dynamics {\eqref{dynsys}} are locally input strictly passive {around} the constant input {$-\omega^*_j = 0$} and {the equilibrium state} $x^{s,*}_j$,}
in the sense of Definition~\ref{Passivity_Definition}, {with {storage function $V^S_j(x^s_j)$} having a strict minimum at {$x^s_j = x^{s,*}_j$}.}
\end{assumption}

\begin{remark}
Assumption \ref{assum_passivity} is a decentralized condition on the continuous power supply dynamics
{that holds} for many dynamical systems considered in the literature and allows the inclusion of higher order schemes (see also \cite{kasis2016primary} and the discussion within it). In Section \ref{Sec: Discussion}, we provide various relevant examples that satisfy Assumption \ref{assum_passivity} to demonstrate its practicality within the considered secondary frequency control setting.
\end{remark}

\section{{{On-Off Loads}}}\label{Sec:switch}
\subsection{Problem formulation}
\label{Sec:switch_Problem}

In this section, we consider {frequency-dependent} on-off loads that respond to frequency deviations by switching to an appropriate state in order to aid the network at urgencies. As the network returns to its {normal operating conditions}, the loads return to their initial state
 as well, hence affecting users comfort for short periods only. {In particular, for {each}
$j \in N$, we {let {the} controllable demand {deviations} {$d^c_j$} be {given}
by the discontinuous map $f^c_j:\mathbb{R} \rightarrow \mathbb{R}$ {defined as}}}
\begin{align}\label{sys2_dc}
{d^c_j= f^c_j(\omega_j)} &=
\begin{cases}
\overline{d}_j, \quad \omega_j > \overline{\omega}_j,\\
\;0, \quad \underline{\omega}_j < \omega_j \leq \overline{\omega}_j,\\
\underline{d}_j, \quad \omega_j \leq \underline{\omega}_j,
\end{cases}
\end{align}
where  $-\infty < \underline{d}_j\leq 0 \leq \overline{d}_j < +\infty $, and  $\overline{\omega}_j > 0 > \underline{\omega}_j$. The {{map} \eqref{sys2_dc} is} depicted on Figure \ref{figure_switch}. Note that {\eqref{sys2_dc}} can be trivially extended to include more discrete values, that would possibly respond to higher frequency deviations {and could correspond to various controllable loads within a bus with dissimilar switching thresholds.
The extension of our results to such cases is trivial and is omitted for simplicity.}

\begin{figure}[t]
\centering
\includegraphics[trim = 0mm 0mm 0mm 0mm, height = 2.00in,width=2.9in,clip=true]{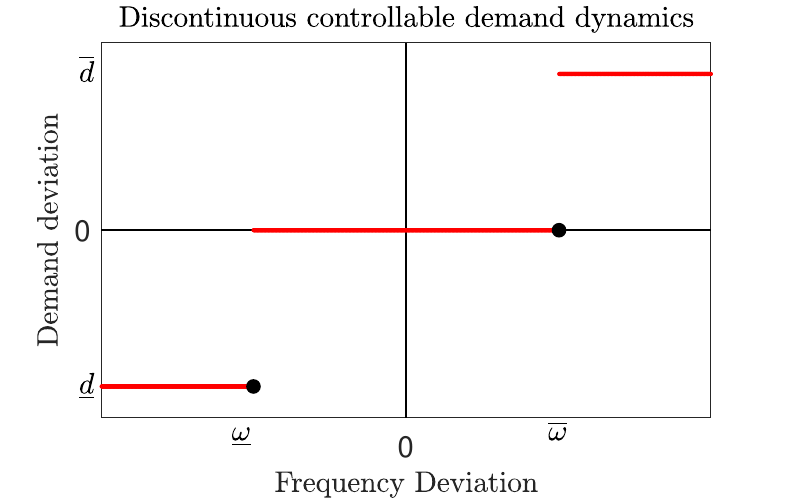}
\caption{{On-off controllable} {demand deviations} as described by \eqref{sys2_dc}.}\label{figure_switch}
\vspace{-0mm}
\end{figure}

{To cope with the {discontinuous nature} of the loads, and to have well-defined solutions to system {\eqref{sys1}--\eqref{dynsys}, \eqref{sys2_dc}} for all time,
{a common approach
is to relax \eqref{sys2_dc} at points of discontinuity with a corresponding map, known as the Filippov set valued map \cite{cortes2008discontinuous}. In particular, this takes here the form}}
\begin{equation}\label{sys2_dc_set_valued}
F[d^c_j]  =  \begin{cases}
[0, \overline{d}_j],\quad \omega_j = \overline{\omega}_j \\[1mm]
[\underline{d}_j,0], \quad\omega_j = \underline{\omega}_j, \\[1mm]
\{{f^c_j(\omega_j)}\}, \, {\rm otherwise.}
\end{cases} j\in N.
\end{equation}

{The state of the} interconnected system \eqref{sys1}--\eqref{dynsys},\eqref{sys2_dc} is denoted by ${x = (\eta, \omega, x^{s})}\in \mathbb{R}^n$, $n = |E| + |N| + \sum_{j \in N} n_j$, where any variable without subscript represents a vector with all respective components.
{By replacing $d^c_j$ in \eqref{sys1} with \eqref{sys2_dc_set_valued} we obtain the following system representation}
\begin{equation}\label{sys_Filippov_representation}
\dot{x} \in Q(x)
\end{equation}
where {$Q$ is a set valued map given by
\vspace{-2mm}
\[
Q(x)=\begin{cases}
 \{\omega_i - \omega_j\}, \; (i,j) \in E, \\[1mm]
\left\{\frac{1}{M_j} (- p_j ^L + s_j -v_j  - \sum_{k:j\rightarrow k} p_{jk}  \right. \\ \left. + \sum_{i:i\rightarrow j} p_{ij}):
\;\; v_j\in F[d^c_j]\right\}, j \in N,
 \\[1mm]
 \{ f_j(x^s_j, -\omega_j) \}, \; j \in N.
\end{cases}
\]}
{
\begin{remark}\label{r:Filip}
{For the analysis of system \eqref{sys1}--\eqref{dynsys}, \eqref{sys2_dc} we will be considering its Filippov solutions. In particular, a Filippov solution of \eqref{sys1}--\eqref{dynsys}, \eqref{sys2_dc} on {an interval $[0,t_1]$ is an absolutely continuous map $x(t)$,  $x: [0,t_1]\rightarrow\mathbb{R}^n$
that satisfies \eqref{sys_Filippov_representation} for almost all $t\in[0,t_1]$.}}
These are often used for the analysis of  {discontinuous systems} as they allow to
overcome the complications associated with the discontinuity\footnote{{The word switching is sometimes used within the paper to refer to the transition of the loads between the on and off states. It should be noted, though, that system \eqref{sys1}--\eqref{dynsys}, \eqref{sys2_dc} is not a switching system as described by e.g. \cite{liberzon2003switching}, but rather a continuous system with a discontinuous right hand side.}} of the vector field.
{Krasovskii solutions \cite[Dfn. 4.2]{goebel2012hybrid} are also frequently used for the study of {discontinuous systems}.
{It should be noted that for the system} \eqref{sys1}--\eqref{dynsys}, \eqref{sys2_dc}, Filippov and Krasovskii solutions are equivalent\footnote{{That is, a Krasovskii solution is also a Filippov solution, and a Filippov solution is also a Krasovskii solution}.} \cite{hajek1979discontinuous}.}
\end{remark}
}

\subsection{{Equilibria and existence of solutions}}

{First, we investigate the existence of Filippov solutions of the system {\eqref{sys1}-\eqref{dynsys}, \eqref{sys2_dc}}, and then study its equilibria.
Existence of Filippov solutions is established by the following lemma.}

\begin{lemma}\label{Uniqueness_Existence}
There exist Filippov solutions of {the system \eqref{sys1}-\eqref{dynsys}, \eqref{sys2_dc}}
{starting from any initial condition $x_0 \in \mathbb{R}^n$.}
\end{lemma}

{\emph{Proof:}}
{The Lemma follows from Proposition 3 in \cite{cortes2008discontinuous}. This states that solutions exist if $Q$ is locally essentially bounded. The latter follows from the boundedness of the step size of the discontinuities in \eqref{sys2_dc} and the Lipschitz property of the rest {of the} dynamics. \hfill $\blacksquare$
}

{An equilibrium of \eqref{sys_Filippov_representation} is defined as follows:}

\begin{definition} \label{eqbrdef}
The {point} $x^* = (\eta^*, \omega^*, x^{s,*})$ defines an equilibrium of the system~\eqref{sys_Filippov_representation} if {${\vect{0}_n} \in Q(x^*)$}.
\end{definition}

{A Lyapunov stable equilibrium point is defined as follows:
\begin{definition}\label{Lyap_stable}
{An equilibrium point $x^\ast$ of \eqref{sys_Filippov_representation} is Lyapunov stable if for all $\epsilon>0$ there exists a $\delta>0$ s.t. any Filippov solution $x(t)$ of  \eqref{sys1}-\eqref{dynsys}, \eqref{sys2_dc} with initial condition $x(0)=x_0$, $\|x_0-x^\ast\|<\delta$, satisfies $\|x(t)-x^\ast\|<\epsilon$ for all $t\geq0$.}
\end{definition}
}

{At an equilibrium of the system, the controllable demand takes its value from a set that depends on $\omega_j^*$, i.e.,} $d^{c,*}_j \in F[d^c_j](\omega_j^*), j \in N$.
{Lemma \ref{eqbr_lemma} below  shows that} {this set is singleton, namely  $Q(x^\ast)=\{{\vect{0}}_n\}$, and $\omega^* = \vect{0}_{|N|}=d^{c,*}$.}

\begin{lemma}\label{eqbr_lemma}
{Let Assumptions  {\ref{assum_eqbl_gen} and \ref{assum_nominal_freq}
hold.}} Then, {there exists an equilibrium point $x^* = (\eta^*, \omega^*, x^{s,*})$ of \eqref{sys_Filippov_representation}.}
 Furthermore, for any equilibrium point of \eqref{sys_Filippov_representation}, {we have}
$\omega^* = \vect{0}_{|N|}$
{and $Q(x^\ast)={\{\vect{0}_n\}}$.}
\end{lemma}
{
\emph{Proof:}}
{
When Assumptions \ref{assum_eqbl_gen} and \ref{assum_nominal_freq} hold, then there exists an equilibrium to system \eqref{sys1}--\eqref{dynsys} with $d^{c}(t) = \vect{0}_{|N|}, t \in \mathbb{R}$ that satisfies $\omega^* = \vect{0}_{|N|}$.
From \eqref{sys2_dc}, $\omega^* = \vect{0}_{|N|}$ implies $d^{c,*} = \vect{0}_{|N|}$, and hence the equilibrium of \eqref{sys1}--\eqref{dynsys} with $d^{c}(t) = \vect{0}_{|N|}, t \in \mathbb{R}$  is also an equilibrium to \eqref{sys_Filippov_representation}.
From  equation \eqref{sys1a} and Assumption \ref{assum_nominal_freq} it also follows that any equilibrium of \eqref{sys_Filippov_representation} satisfies $\omega^* = \vect{0}_{|N|}$ which implies  $Q(x^\ast)={\{\vect{0}_n\}}$ from \eqref{sys2_dc_set_valued}.
\hfill $\blacksquare$
}

\subsection{{Convergence analysis}}

{We {now state the main result of this section. {Its proof is provided in the} {Appendix.}}}

\begin{theorem} \label{convthm}
{{Let {Assumptions~\ref{assum_eqbl_gen}, \ref{assum_0},~\ref{assum_nominal_freq}} hold and consider an equilibrium $x^*$ where {Assumptions \ref{assum1},~\ref{assum_passivity}} hold.}}
Then, {$x^{\ast}$ is Lyapunov stable. Furthermore,} there exists {a compact} neighborhood  {of}
this equilibrium such that Filippov solutions {$x=(\eta, \omega, {x^s}), x:[0,\infty)\rightarrow\mathbb{R}^n$}
of
{the system \eqref{sys1}-\eqref{dynsys}, \eqref{sys2_dc}} starting in {this neighborhood}
{converge} to {the set of equilibria of the system that lie in this neighbourhood}. In particular, the frequency vector $\omega$ converges to $\omega^\ast=\vect{0}_{|N|}$.
\end{theorem}

{Theorem \ref{convthm} states that solutions that start sufficiently close to $x^\ast$ are guaranteed to converge to {the set of equilibrium points},
{which suggests by Lemma \ref{eqbr_lemma} that {the}} frequency is restored to its nominal value after {transient} load-side participation.}
{{Note} that Theorem \ref{convthm} concludes convergence {to a set of equilibria rather than a single equilibrium point (see also Remark \ref{Rem:semistability} after the proof of the theorem)}.
}
 {It should also be {noted} that the local nature of the {convergence result} in Theorem \ref{convthm} {is due to the
{nonlinearity of the}} sinusoids in the {power transfer} \eqref{sys1d} and the fact that {the passivity assumption on the generation dynamics (Assumption~\ref{assum_passivity}) holds locally\footnote{{In particular, if the sinusoids where linearized and the passivity property on the generation dynamics was global then a set of equilibrium points of the resulting dynamical system would have been globally attractive.}}}}.

\subsection{{Chattering}}\label{sec: Zeno}

{A possible feature of {discontinuous systems}}
is the occurrence of infinitely many {transitions} within some finite time, a phenomenon known as {chattering}, e.g. \cite{lee2007chattering}.
Such behavior is {often undesirable and impedes practical implementations.}

{In our setting, {chattering} may occur in on-off loads} as shown numerically in Section \ref{Simulation}.
The reason such behavior may occur is that
the frequency derivative may change sign when passing a discontinuity, causing the vector field to point towards the discontinuity.
{For instance, suppose that $0< {M_j}\dot{\omega}_j(t_1)<\overline{d}_j$ for some time instant $t_1>0$, and that the threshold $\overline{\omega}_j$ is met at this time. Then, the load $d^c_j$ switches on, causing a sign change in the value of $\dot{\omega}_j$.
{Hence, the frequency vector field will point at a direction of frequency decrease that will force the load to switch off.
These on/off switches} occur infinitely many times in a finite time, resulting in the aforementioned {chattering} behavior. Note that this phenomenon is only observed here during the transient response of the loads, as the mechanical power injection will eventually {dictate the sign of the vector field and regulate} the frequency to its nominal value as shown in Theorem \ref{convthm}.}

\section{Hysteresis on controllable loads}\label{Sec:Hysteresis}

{
\subsection{Problem formulation}\label{sec_problem_formulation}

In this section, we propose the use of hysteretic dynamics in on-off controllable loads, which means that a controllable load switches on and off at different frequency thresholds. As {it will be shown}, this modification will ensure that the system does not exhibit {chattering}.
For relevant applications of hysteric dynamics in ruling out {chattering and other {undesired} features}, see e.g. \cite{lee2007chattering,mayhew2011quaternion,ceragioli2011discontinuities,prieur2005asymptotic}.

We consider {hysteretic dynamics for controllable loads that {satisfy}}
\begin{subequations}\label{sys_hysteresis}
\begin{equation}
d^c_j = \overline{d}_j \sigma_{j}
\end{equation}
\begin{equation}\label{sys_hysteresis_g}
\begin{aligned}
\sigma_{j}(t^+) { = } {\begin{cases}
{\rm sgn}(\omega_j), \hspace{2mm} \qquad |\omega_j| > \omega^1_j\\[1mm]
0,  \hspace{19mm} |\omega_j| < \omega^0_j\\[1mm]
\sigma_j(t),  \hspace{14mm} \omega^0_j < |\omega_j| < \omega^1_j
\end{cases}}
\end{aligned}
\end{equation}
\end{subequations}
where $j\in N$, $t^+ = \lim_{\epsilon \rightarrow 0} (t + \epsilon)$, and the frequency thresholds $\omega_j^0, \omega_j^1$, satisfy
$\omega_j^1 >\omega_j^0>0$.
 {Note that $\sigma_j$ takes its value from the set $P=\{-1,0,1\}$. A schematic of the hysteretic dynamics is shown in Figure \ref{hysteresis_figure}.

 In the remainder of this section and in the next section we provide a detailed definition of the overall dynamical system as a hybrid system \cite{goebel2012hybrid} and use corresponding tools for its analysis.}
}

\begin{figure}[t]
\centering
\includegraphics[trim = 0mm 0mm 0mm 0mm, height = 2.00in,width=2.9in,clip=true]{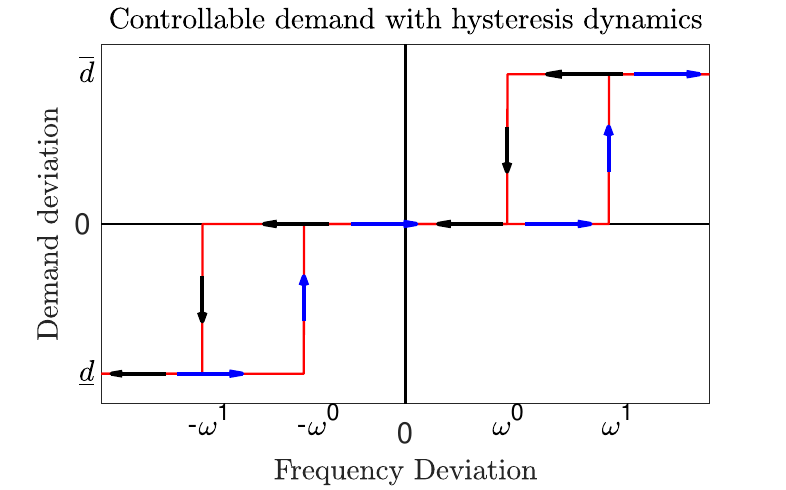}
\caption{Hysteresis dynamics for controllable loads as described by \eqref{sys_hysteresis}.}
\label{hysteresis_figure}
\end{figure}

{In particular, system} {\eqref{sys1}--\eqref{dynsys},} {with hysteretic controllable loads as in \eqref{sys_hysteresis},} can be described by the states $ z = (x,\sigma$), where $x = (\eta, \omega, x^{s}) \in \mathbb{R}^n$ is the continuous state, and $\sigma \in P^{|N|}$ the discrete state. {We also denote by $\Lambda$ the domain where the state $z$ takes values, which is a subset of $\mathbb{R}^n \times P^{|N|}$, where for each} {$\omega_j$ the value of $\sigma_j$ satisfies   $\sigma_j \in \mathcal{I}_j(\omega_j)$, where
\[
\mathcal{I}_j(\omega_j)= \begin{cases}
\{{\rm sgn}(\omega_j)\}, {\qquad |\omega_j| >  \omega^1_j,}\\[1mm]
\{0\},  \hspace{17mm} {|\omega_j| <  \omega^0_j,}\\[1mm]
\{0, {\rm sgn}(\omega_j)\}, \quad \omega^0_j \leq |\omega_j| \leq \omega^1_j.
\end{cases}
\]
This restriction on the state space ensures that the initial conditions for $\sigma$ are compatible with {their right hand limits in} \eqref{sys_hysteresis_g} and is made to simplify the presentation of the subsequent analysis and results.
}

{The continuous part of the dynamics is given by
{\begin{subequations}\label{sys4_hysteresis}
\begin{equation}\label{sys4a}
\dot{\eta}_{ij} = \omega_i - \omega_j, \; (i,j) \in E,
\vspace{-4mm}
\end{equation}
\begin{multline}\label{sys4b}
M_j \dot{\omega}_j = - p_j ^L + s_j - (\overline d_j\sigma_j)\\
 - \sum_{k:j\rightarrow k} p_{jk} + \sum_{i:i\rightarrow j} p_{ij}, \; j\in N,
\end{multline}
\vspace{-2mm}
\begin{equation}\label{sys4d}
\dot{x}^s_j = f_j(x^s_j,-\omega_j), \; j \in N,
\end{equation}
\begin{equation}\label{sys4e}
\dot{\sigma}_j = 0, j \in N,
\end{equation}
\end{subequations}
{where $p_{ij}$ and $s_j$ are given by \eqref{sys1d} and \eqref{dynsys}, respectively.}}
{This is} valid when $z$ belongs to the set
\begin{equation}\label{e:C}
{C=\Lambda}=\{ z \in \mathbb{R}^n \times P^{|N|}:  \sigma_j \in \mathcal{I}_j(\omega_j), \;\forall j\in N\}.
\end{equation}

Alternatively, when {$z$ belongs to the set $D$ defined {as
\begin{equation}\label{set_D}
D = \{ z \in \Lambda: |\omega_j| \in \{\omega^0_j, \omega^1_j\},  \sigma_j \in \mathcal{I}^D_j(\omega_j), \;\forall j\in N\}
\end{equation}
{with}}
\[
\mathcal{I}^D_j(\omega_j)= \begin{cases}
\{0\}, \qquad \qquad \hspace{2.5mm} {|\omega_j| =  \omega^1_j},\\[1mm]
\{{\rm sgn}(\omega_j)\},  \hspace{7mm} {|\omega_j| =  {\omega^0_j},}\\[1mm]
\end{cases}
\]
the} system dynamics evolve according to the following discrete update rule:

\begin{subequations}\label{sys4_g}
\begin{equation}
 x^+ = x, \quad \\
\end{equation}
\begin{equation}
\sigma_{j}(t^+) = \begin{cases}
{\rm sgn}(\omega_j), \qquad {|\omega_j| =  \omega^1_j},\\[1mm]
0,  \hspace{17mm} {|\omega_j| =  \omega^0_j}.\\[1mm]
\end{cases}
\end{equation}
\end{subequations}

 We can now provide the following compact representation for the hybrid system {\eqref{sys1}--\eqref{dynsys},~\eqref{sys_hysteresis},}
\begin{subequations}\label{sys4}
\begin{equation}\label{e:f(z)}
{\dot{z} = f(z), \quad z \in C,}
\end{equation}
\begin{equation}
{z^+ = g(z), \quad z \in D,}
\end{equation}
\end{subequations}
where {the maps $f(z): C \rightarrow \Lambda$ and $g(z): D \rightarrow C$  are given by {\eqref{sys4_hysteresis}} and \eqref{sys4_g}, respectively, the sets $C$ {and $D$ by \eqref{e:C} and \eqref{set_D} respectively,
 and $z^+=z(t^+)$.}}
Note that $z^+ = g(z)$ represents a discrete dynamical system where $z^+$ is determined by the current value of the state $z$ and the update rule given {by $g$.}

\subsection{Analysis of equilibria and solutions}

{Before investigating stability of the hybrid system in \eqref{sys4}, we {characterize} its equilibria, and establish existence and completeness of solutions.}

{Note that we call a point $z^* = (x^*, \sigma^*)$ an equilibrium of \eqref{sys4} {if $f(z^*)=0$ {when $z^* \in C$, and $z^\ast=g(z^\ast)$ when $z^* \in D$.}
Now, we state the following lemma:

\begin{lemma}\label{eqlbr_hysteresis}
{Let {Assumptions}  {\ref{assum_eqbl_gen} and \ref{assum_nominal_freq}
hold.}} {Then, there exists an equilibrium of \eqref{sys4}.
Furthermore, }
 for any equilibrium point $z^* = (x^*, \sigma^*)$ of \eqref{sys4}, we have $\omega^* = \sigma^* = \vect{0}_{|N|}$. Moreover, $z^* \in C$.
\end{lemma}
{\emph{Proof:}}
{The existence of an equilibrium of \eqref{sys4} follows analogously with the proof of Lemma \ref{eqbr_lemma}.}
{Recall that any equilibrium $z^*$ of \eqref{sys4} satisfies} $f(z^*) = 0$ {when $z^* \in C$, and $z^\ast=g(z^\ast)$ when $z^* \in D$.}
{Since $g(z):D \rightarrow C$, it therefore follows that $z^* \in C$.}
{{From Assumption~\ref{assum_nominal_freq}}}
and equations \eqref{sys1a} and \eqref{sys1d} at equilibrium, it follows that $\omega^* = \vect{0}_{|N|}$, which implies that $\sigma^* = \vect{0}_{|N|}$.
\hfill $\blacksquare$
}}

{{Following \cite{goebel2012hybrid}}, we provide a definition of a hybrid time domain, hybrid solution and complete and maximal solutions for systems described by \eqref{sys4}.

\begin{definition}{(\cite{goebel2012hybrid})}\label{dfn_hybrid_domain_solution}
{A subset of {$\mathbb{R}_{\geq0} \times \mathbb{N}_0$} is a hybrid time domain if it is a union of a finite or infinite sequence of intervals  $[t_\ell, t_{\ell+1}] \times \{\ell\}$, with the last interval (if existent) possibly of the form {$[t_\ell, t_{\ell+1}] \times \{\ell\}$, $[t_\ell, t_{\ell+1}) \times \{\ell\}$, or $[t_\ell, \infty) \times \{\ell\}$.}
Consider a function $z(t,\ell): K \rightarrow \mathbb{R}^n$ defined on a hybrid time domain $K$ such that for every fixed $\ell \in \mathbb{N}$, $t\rightarrow z(t,\ell)$ is locally absolutely continuous on the interval $T_\ell = \{t: (t,\ell) \in K\}$.
 The function $z(t,\ell)$ is a solution to the hybrid system $\mathcal{H} = (C,f,D,g)$ if {$z(0,0) \in {{C} \cup D}$}, and for all $\ell \in \mathbb{N}$ such that $T_\ell$ has non-empty interior {(denoted by ${\rm int}T_l$)}
\begin{align*}
& z(t,\ell) \in C, \text{ for all t} \in {{\rm int}T_l}, \\
& \dot{z}(t,\ell) \in f(z(t,\ell)), \text{ for almost all } t \in T_\ell,\\
& \hspace{-12mm} \text{and for all } (t,\ell) \in K \text{ such that } (t,\ell+1) \in K,\\
& z(t,\ell) \in D, \; z(t,\ell+1) \in g(z(t,\ell)).
\end{align*}
A solution $z(t,\ell)$ is complete if $K$ is unbounded. {A solution {$z$ is maximal} if there does not exist another {solution $\tilde z$} with time domain $\tilde K$ such that $K$ is a {proper} subset of $\tilde K$ and $z(t,j) = \tilde{z}(t,j)$ for all {$(t,j) \in K$}}.}
\end{definition}

For convenience in the presentation we will refer to maximal solutions by just solutions.}}
Existence of complete solutions to \eqref{sys4} are {established in} the following lemma.
\begin{lemma}\label{lemma_existence_uniq_hysteresis}
{There exists a complete solution $z=(x, \sigma)$ to \eqref{sys4}, starting from any initial condition $z(0,0) \in {\Lambda}$.}
\end{lemma}

Furthermore, the following {proposition} shows the existence of some finite dwell time between switches of states $\sigma_j$ {for any bounded solution.  Within it, we denote  the time-instants where the value of $\sigma_j$ changes by $t_{\ell,j}, \ell \in {\mathbb{N}_0,} j \in N$.}

\begin{proposition}\label{dwell_time_lemma}
For any complete bounded solution of \eqref{sys4}, there {exists $\tau > 0$}  such that {$\min_{{\ell} \geq 1} (t_{{\ell}+1, j} - t_{{\ell}, j}) \geq \tau$} for any $j \in N$.
\end{proposition}

\begin{remark}
The importance of Proposition \ref{dwell_time_lemma} is that it shows that no {chattering} will occur for any complete bounded solution of system \eqref{sys4}. This is because for any finite time interval {$\tau$}, the vector $\sigma$ changes at most $|N|$ times. This highlights the practical advantage of \eqref{sys4} compared to \eqref{sys_Filippov_representation}. This analytic result is verified by numerical simulations in a realistic power network, as discussed in section \ref{Simulation}.
\end{remark}
\vspace{-0mm}
\subsection{Stability of hysteresis system}\label{sec: Stability_hybrid}

{Now, we are at the position to state the stability result concerning the system \eqref{sys4}.}
\begin{theorem}\label{conv_thm_hysteresis}
Let {{Assumptions~\ref{assum_eqbl_gen},~\ref{assum_0},~\ref{assum_nominal_freq},} hold and consider an equilibrium {$z^* = (x^*, \sigma^*)$ of \eqref{sys4},
for which} {Assumptions~\ref{assum1},~\ref{assum_passivity}}  hold.}
Then there exists {a compact} {neighborhood of $x^\ast$} such that solutions $z=(x, \sigma)$, with {$x(0)$ in this neighbourhood} and $\sigma(0)\in {\mathcal{I}}(\omega(0))$,  asymptotically converge to {the set of equilibria of  \eqref{sys4} whose continuous state lies in this neighbourhood}. In particular, the vectors $\omega$ and $\sigma$ converge to the vector
$\vect{0}_{|N|}$.
\end{theorem}

\begin{remark}
Theorem \ref{conv_thm_hysteresis} shows that the hysteretic dynamics in \eqref{sys_hysteresis} do not compromise the stability of the system.
{This, together with the absence of {chattering} shown in Proposition \ref{dwell_time_lemma}, promotes the use of hysteretic dynamics as a {means} to provide practical and non-disruptive on-off load side control to the power network.}
\end{remark}

{
\begin{remark}
Although the controllable loads {are at their nominal value} at {equilibrium}, they provide ancillary services to the network, and improve the performance in transient time. This is numerically investigated in Section \ref{Simulation}. To clarify, note that the convergence region in Theorem \ref{conv_thm_hysteresis} is not restricted by the switches, but rather by the nonlinearity of the frequency dynamics.
\end{remark}}

{
\begin{remark}\label{rem: Graph}
We note that no assumption is made on the graph $(G,E)$ except that it is connected. Therefore, all results presented in the current and previous section are applicable on arbitrary connected graph topologies.
\end{remark}
}

\section{Discussion}\label{Sec: Discussion}

{It was shown in section \ref{Sec:Hysteresis} that hysteretic on-off loads can be incorporated in an arbitrary network without compromising stability if the continuous, and potentially higher order, supply dynamics satisfy the conditions described in section \ref{Sec: Passive_dynamics}. In this section we demonstrate the generality and applicability of these conditions by discussing various examples that fit within the proposed framework.}

An important feature of our analysis is that it allows to consider {a} general class of dynamics
where {the} power supply dynamics are described by the summation of a {damping term}
and the {series} interconnection of {a PI controller} {with} an asymptotically stable {general linear system, potentially of higher order.} {The transfer function from $-\omega_j$ to $s_j$ of such a system is given {by} ${S}_j(s) =  (\frac{K_j}{s}+{\tilde{K}_j})G_j(s) + D_j$ {where $K_j, \tilde{K}, D_j$ are positive constants {and\footnote{{We assume that no zero/pole cancellation occurs between $G_j(s)$ and the PI controller.}}} $G_j(s)$ is a transfer function with no unstable poles}.}
{Assumption~\ref{assum_nominal_freq} and the assumptions on \eqref{dynsys}, are clearly satisfied for this system. Moreover,
Assumption~\ref{assum_passivity} is satisfied if there exists $\epsilon_j > 0$ such that the {perturbed} transfer function  $S_j(s)-\epsilon_j$
is positive real.} This can be numerically verified with appropriate LMI (Linear Matrix Inequality) conditions \cite[Lemma 6.2]{khalil1996nonlinear}, {i.e. a computationally efficient convex feasibility problem}, or graphically by examining the Nyquist plot of $\color{black}{S}_j(j\omega)$, which needs to lie within the open right half plane.

First, we consider power supply dynamics where the output response {follows
{from a PI control scheme
acting on a frequency input with a lag}}. Such dynamics can be described~by
\begin{subequations}\label{int_first_order}
\begin{equation}
\dot{\alpha}_j = -K_j \omega_j, \quad
{\tau_{\beta, j} \dot{\beta}_j = -\beta_j + \alpha_j  - {\tilde{K}_j} \omega_j},
\end{equation}
\begin{equation}
s_j = \beta_j - D_j\omega_j,
\end{equation}
\end{subequations}
where $\alpha_j, \beta_j$ are the states of the power supply variables and $K_j, {\tilde{K}_j,} \tau_{\beta_j}, D_j$ are positive constants describing {the} integrator {and droop} gains, the time constant {associated with power generation and the} frequency damping respectively. For \eqref{int_first_order}, it can be shown that Assumption \ref{assum_passivity} is satisfied if {$K_j \tau_{\beta_j} < D_j + {\tilde{K}_j}$}.  Furthermore, Assumption \ref{assum_nominal_freq} is trivially satisfied.

An important aspect of our analysis lies in its ability to consider higher order schemes. A significant example of this can be seen in second order generation dynamics, often considered to describe turbine-governor behaviour (e.g. \cite{Bergen_Vittal}).
Below, we consider the series interconnection of an integrator and a second order system {with} frequency damping, described by
\vspace{-0.2cm}
\begin{subequations}\label{int_second_order}
\begin{equation}
\dot{\alpha}_j = -K_j \omega_j, \quad
\tau_{\beta, j} \dot{\beta}_j = -\beta_j + \alpha_j,
\end{equation}
\begin{equation}
\tau_{\gamma, j} \dot{\gamma}_j = -\gamma_j + \beta_j, \quad
s_j = \gamma_j - D_j\omega_j,
\end{equation}
\end{subequations}
where $\alpha_j, \beta_j, \gamma_j$ are internal states of the system and $K_j, \tau_{\beta_j}, \tau_{\gamma_j}, D_j$ positive constants.  It is straightforward to deduce that such system satisfies Assumption \ref{assum_nominal_freq}. Furthermore, it can be shown that \eqref{int_second_order} satisfies Assumption \ref{assum_passivity} {when $K_j (\tau_{\beta, j} + \tau_{\gamma, j}) < D_j$.}

The Nyquist plots of the dynamics in \eqref{int_first_order}--\eqref{int_second_order} are illustrated {in} Figure \ref{Nyquist_passivity} that depicts the plots of these schemes when the suggested gain conditions are marginally satisfied.
{Figure~\ref{Nyquist_passivity} demonstrates the passivity of the considered {dynamics when the gain conditions are satisfied,} since the Nyquist plots are {in} the right half plane.}

\begin{figure}
\centering
\includegraphics[trim = 5mm 1mm 12mm 5mm, height = 2.00in,width=2.5in,clip=true]{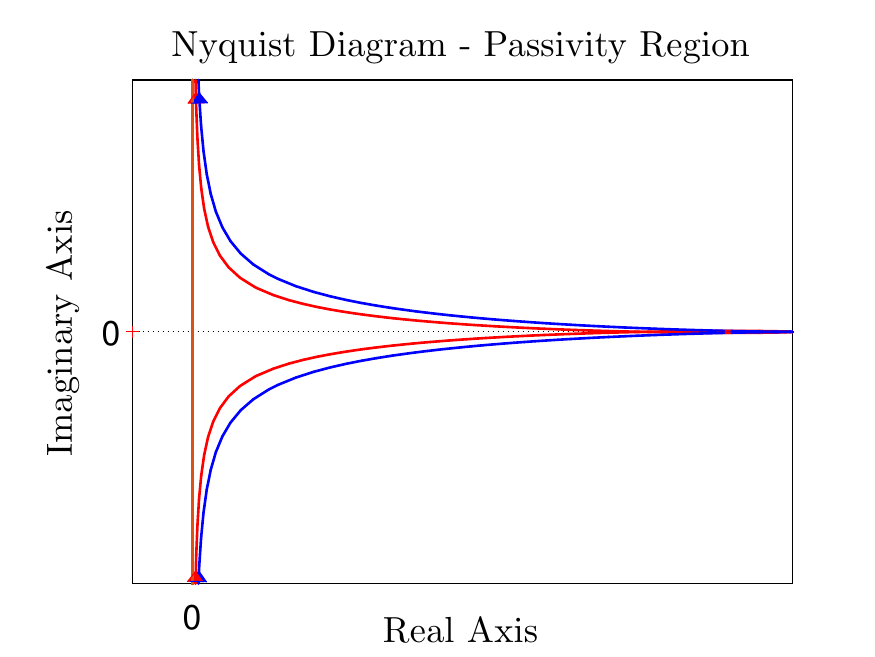}
\vspace{-0mm}
\caption{Nyquist plots of the {systems from $\omega_j$ to $s_j$ described} in \eqref{int_first_order} and \eqref{int_second_order}, illustrated by blue and red lines respectively, when the gain conditions {stated} are marginally satisfied. The passivity of both schemes is implied by the fact that both plots lie on the right half plane.}
\vspace{-0mm}
\label{Nyquist_passivity}
\end{figure}

Furthermore, our framework allows to incorporate various realistic dynamical schemes that satisfy Assumption \ref{assum_passivity}, but not Assumption \ref{assum_nominal_freq}.
{This is relevant since Assumption \ref{assum_nominal_freq} {only} needs to be satisfied by at least one bus in the network, as mentioned in Remark \ref{rem_assum2}.}
To demonstrate an important example of such {a case}, we consider the fifth order dynamics used by the Power System Toolbox~\cite{cheung2009power} to describe turbine governor behavior within the NPCC network. In this model, power supply $\hat{s}_j$ is related to the negative frequency deviation $-\hat{\omega}_j$ {in the Laplace domain} via the transfer function
\[
G_j(s)=K_j\frac{1}{(1+sT_{s,j})}\frac{(1+sT_{3,j})}{(1+sT_{c,j})}\frac{(1+sT_{4,j})}{(1+sT_{5,j})} + D_j,
\]
where $T_{s,j}, T_{3,j}, T_{c,j}, T_{4,j}, T_{5,j}$ are time-constants and $K_j$ and $D_j$ are the droop and damping coefficients respectively. Realistic values for these systems are provided by the Power System Toolbox and it can be shown that Assumption \ref{assum_passivity} is satisfied by $20$ out of the $22$ buses with turbine governor dynamics. Furthermore, for the remaining $2$ buses, Assumption \ref{assum_passivity} is satisfied if the damping coefficients are increased by $28\%$ and $37\%$ respectively. This demonstrates that the passivity conditions are satisfied by existing implementations.
Note that further examples of schemes that satisfy Assumption \ref{assum_passivity}, including static nonlinearities, first and second order systems are provided in \cite{kasis2016primary}.

\section{Simulation on the NPCC 140-bus system} \label{Simulation}

In this section we verify our analytic results with a numerical simulation on the Northeast Power Coordinating
Council (NPCC) 140-bus interconnection system, using the Power System Toolbox~\cite{cheung2009power}. This model is more detailed and realistic than our analytical one, including line resistances, a DC12 exciter model{, a subtransient reactance generator model, and turbine governor dynamics\footnote{The details of the simulation models can be found in the Power System Toolbox data file datanp48.}.

The test system consists of 93 load buses serving different types of loads including constant active and reactive loads and 47 generation buses. The overall system has a total real power of 28.55GW. For our simulation, we added three {uncontrollable} loads on units 8, 9 and 17, each having a step increase of magnitude $2$ p.u. (base 100MVA) at $t=1$ second.

Controllable demand was considered within the simulations, with {controllable} loads controlled every 10ms. Secondary control was performed from the generation side, where the dynamics satisfied the conditions presented in Section \ref{Sec: Passive_dynamics}.

The system was tested at two different cases. In case (i), on-off controllable loads with dynamics as in \eqref{sys2_dc} were included on 20 load buses. The values for $\overline{\omega}_j$ were selected from a uniform distribution within the range $[0.02 \; 0.07]$ and those of $\underline{\omega}_j$ by $\underline{\omega}_j = -\overline{\omega}_j$. Controllable loads were also included on 20 load buses for case (ii), but with dynamics  described by \eqref{sys_hysteresis}. To have a fair comparison, the same frequency thresholds were used for both cases, with $\omega^1_j = \overline{\omega}_j$ and $\omega^0_j = 0.15\omega^1_j$. Also, $\overline{d} = -\underline{d} = 0.25p.u.$ was used for both cases. Furthermore,  in order to account for any potential effect of the reactive power deviations in the response of the network, {on-ff} loads were assumed to induce a reactive power effect of magnitude $0.25p.u.$, i.e. assumed a power factor of about $0.7$ which is in most cases {an underestimate}.
Cases (i) and (ii) will be referred to as the 'switching' and 'hysteresis' cases respectively.

The frequency at bus 99 for the two tested cases is  shown in Fig. \ref{Frequency}, where it can be seen that the frequency returns to its nominal value for both cases, as suggested in Theorems~\ref{convthm} and \ref{conv_thm_hysteresis}. Moreover, Fig. \ref{overshoot_freq} demonstrates that the inclusion of {on-off} loads decreases the maximum overshoot in frequency, by comparing the largest deviation in frequency with and without on-off controllable loads at buses $1-40$, where frequency overshoot was seen to be the largest.
From Figure \ref{overshoot_freq}, it can be seen that the largest frequency overshoot has dropped from about 0.28$Hz$ to about 0.18$Hz$, demonstrating a drop of more {than} $30\%$. Note that the similar responses from switching and hysteretic loads on Fig. \ref{overshoot_freq} follow from the fact that identical frequency thresholds have been chosen. Furthermore, the simulation results demonstrated that deviations in voltage magnitude  have always been within $3.3\%$ of the nominal value at all times, justifying the {constant voltages} assumption in Section \ref{sec:Network_model}.

  Figure \ref{Fig_switch_case} shows controllable demand at 4 buses for case (i), depicting very fast switches, indicated by the thick blue lines, which demonstrates {chattering} behavior\footnote{Note that in practice loads won't be able to switch arbitrarily fast, since their physical parameters will impose some time delays between switches. However, very fast switching will still be observed, which is an undesirable feature.}. In contrast, when case (ii) is considered, such fast switching in loads is not observed on those 4 buses, as exhibited in Figure \ref{Fig_hysteresis_case}.
  Furthermore, it was seen within the simulations that all loads switched off after $30s$, which demonstrates the non-disruptive nature of the two schemes, since loads return to their nominal demand after a brief period. It should be noted that switching loads demonstrated a faster restoration than hybrid ones as a result of their higher switch off thresholds.
Thus, this numerical investigation supports the analysis of the paper,  verifying that {frequency-dependent} on-off loads are able to provide ancillary service in the power network and that hysteresis eliminates any {chattering} behavior.

\begin{figure}[t!]
\centering
\includegraphics[trim = 0mm 0mm 0mm 0mm, height = 2.20in,width=3.15in,clip=true]{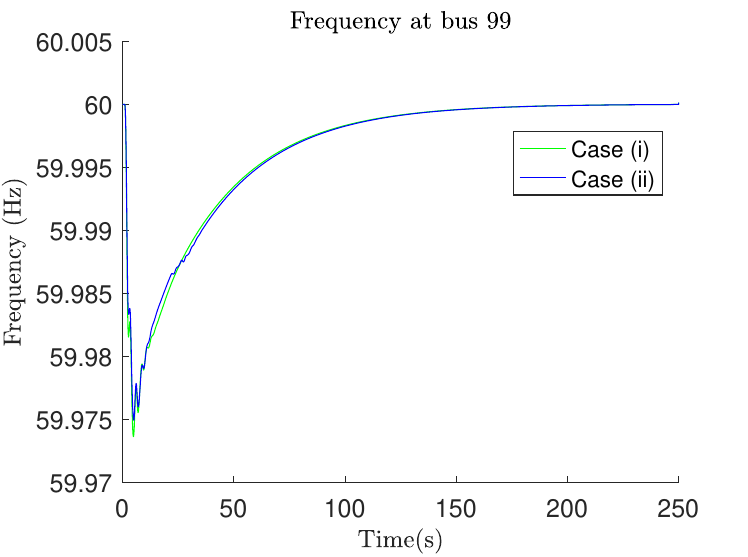}
\caption{Frequency at bus 99 with controllable load dynamics as in the following two cases: i) Switching case, ii) Hysteresis case.}
\label{Frequency}
\vspace{-0mm}
\end{figure}

\begin{figure}[t!]
\centering
\includegraphics[trim = 1mm 0mm 0mm 0mm, height = 2.20in,width=3.15in,clip=true]{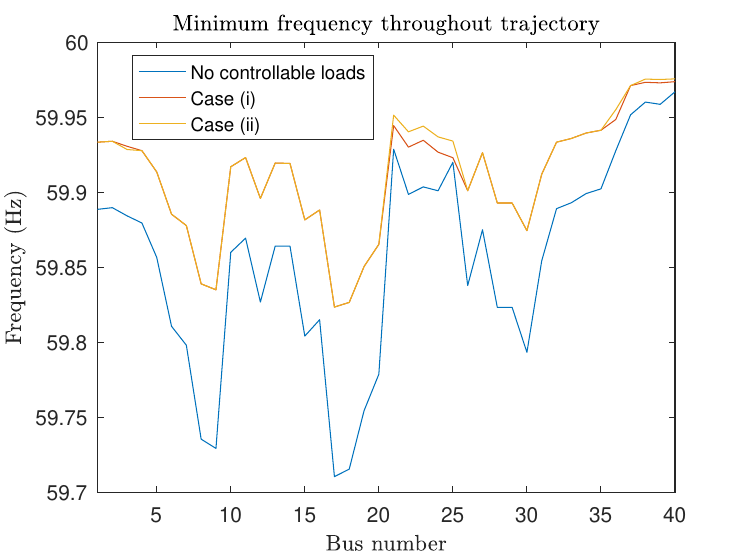}
\caption{Largest frequency overshoot for buses $1-40$ for three cases: (i) Use of switching loads, (ii) Use of hysteresis loads, (iii) No use of controllable loads. {Note that the graphs for cases (i), (ii) are almost identical in the figure.}}
\label{overshoot_freq}
\vspace{-0mm}
\end{figure}

\begin{figure}[t!]
\centering
\includegraphics[trim = 5mm 3mm 10mm 3mm, height = 2.2in,width=3.15in,clip=true]{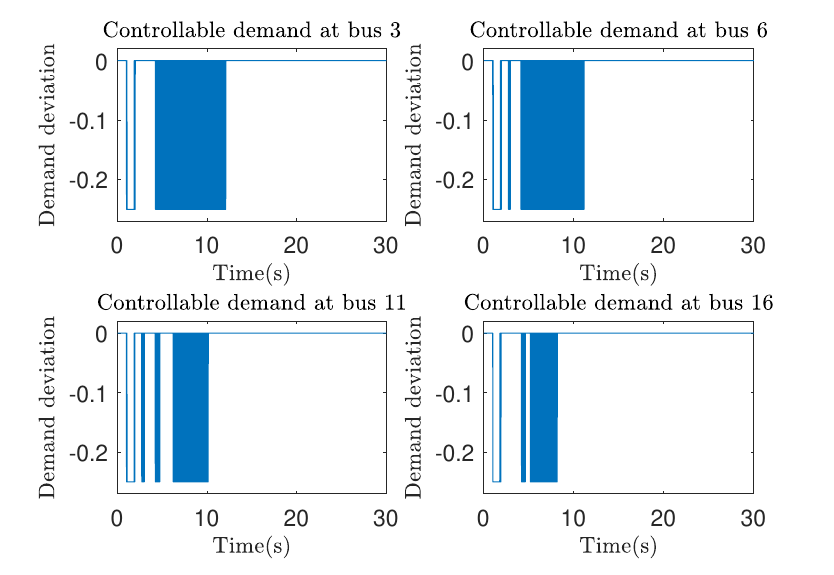}
\caption{Controllable demand at 4 buses with Switching on-off loads.  All loads stay switched off after 15$s$.}
\label{Fig_switch_case}
\vspace{-0mm}
\end{figure}

\begin{figure}[t!]
\centering
\includegraphics[trim = 5mm 3mm 10mm 3 mm, height = 2.2in,width=3.15in,clip=true]{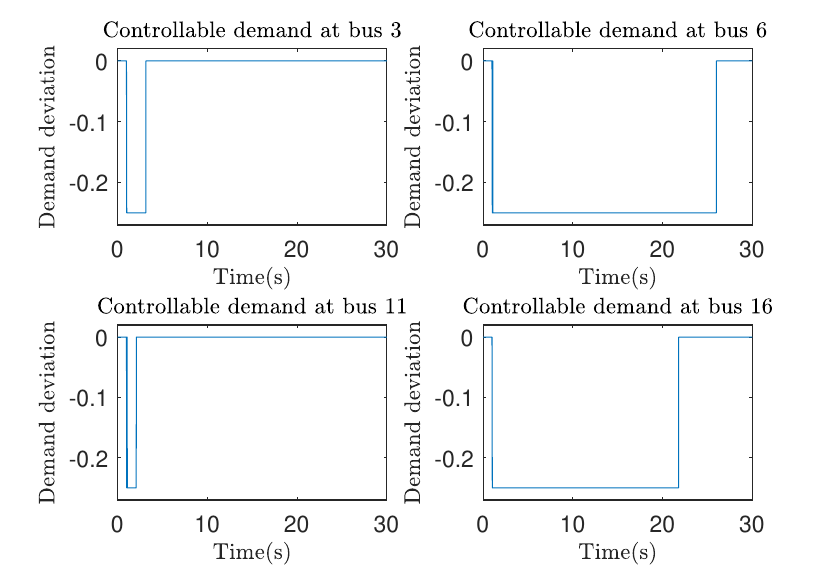}
\caption{Controllable demand at 4 buses with Hysteresis on-off loads. All loads stay switched off after 30$s$.}
\label{Fig_hysteresis_case}
\vspace{-0mm}
\end{figure}

\section{Conclusion}\label{Conclusion}
We have considered the problem of secondary frequency control {in power networks} where controllable on-off loads provide ancillary services.
We first considered a general class of continuous dynamics {together with} loads that switch on when some frequency threshold is reached and off otherwise. Stability guarantees are provided for this framework. Furthermore, it is discussed that such on-off schemes might exhibit arbitrarily fast switching, which might limit their practicality. To cope with this issue, on-off loads with hysteretic dynamics were considered. It has been shown that such loads do not exhibit {chattering} and that their inclusion does not compromise power network stability. Hence, such schemes are usable for practical implementations. Both schemes ensure that controllable loads return to their nominal behavior at equilibrium and hence that disruptions occur for brief periods only.
We discuss that our analysis allows the incorporation of on-off loads {together with} a broad range of continuous dynamics, including various highly relevant and practical schemes which highlights the applicability of our results.
Numerical simulations on the NPCC 140-bus system verify our main findings, demonstrating} that the presence of on-off loads reduces the frequency overshoot and that hysteresis schemes resolve issues caused by {chattering}.
Interesting potential extensions in the analysis include incorporating {frequency-dependent} on-off loads within the primary frequency control timeframe, more advanced load dynamics, and models that take into account voltage dynamics.

\section*{Appendix: proofs of results}\label{appendix}

Within the proofs  of Lemma \ref{eqbr_lemma} and Theorem \ref{convthm} we will make use of the following equilibrium equations for system \eqref{sys1}--\eqref{dynsys}, which follow from Definition \ref{eqbrdef} and Lemma \ref{eqbr_lemma}.  Below, we let {$\bar{N}\subseteq N$} be the set of all buses that satisfy the condition in Assumption \ref{assum_nominal_freq}, i.e.  those that satisfy $f_j(\bar{x}^s_j, -\bar{\omega}_j) = 0$ only when $\bar{\omega}_j = 0$.
\begin{subequations} \label{eqbr}
\begin{gather}
0 = \omega^*_i - \omega^*_j, \; (i,j) \in E, \label{eqbr1} \\
0 = - p_j^L \! + \! s_j^{*} \!  - \!\! \sum_{k:j\rightarrow k} p^*_{jk} + \!\! \sum_{i:i\rightarrow j} p^*_{ij}, \; j\in N, \label{eqbr2} \\
{p^*_{ij}=B_{ij} \sin \eta^*_{ij}}, \; (i,j) \in E,\ \label{eqbr4} \\
0 = \omega^*_j, j \in \bar{N}, \label{eqbr5}\\
 0 = d^{c,*}_j, j \in {N} \label{eqbr6}.
\end{gather}
\end{subequations}

\emph{Proof of Theorem~\ref{convthm}:}
{To prove Therorem~\ref{convthm}, we will make use of \cite[Thm. 3]{bacciotti1999stability}.
First, note that the set valued map $Q$ in \eqref{sys_Filippov_representation} \ak{is upper semicontinuous\footnote{\ak{For the definition of an upper semicontinuous function, see \icl{e.g.} \cite[p.49]{cortes2008discontinuous}.}}
 and} {takes compact,} convex values.
{For}} the dynamics \eqref{sys1}--\eqref{dynsys}, and \eqref{sys2_dc}, we define
\begin{align}\label{V_Lyapunov_2}
&\bar{V}(\eta, \omega) = V_F(\omega) + V_P(\eta)
\end{align}
where $V_F (\omega) = \frac{1}{2}\sum_{j \in N} M_j \omega_j^2$,
and
$
V_P(\eta) = \sum_{(i,j) \in E} B_{ij} \int_{\eta^*_{ij}}^{\eta_{ij}} ( \sin \phi - \sin \eta^*_{ij} ) \, d\phi.
$
}

{By substituting~\eqref{sys1a} and \eqref{sys1b} for $\dot{\eta}_{ij}$ and $\dot{\omega}_j$  and using the differential inclusion for $d^c_j$ for $j \in N$, the {set-valued} time-derivative of $\bar{V}$ along {the solutions of~\eqref{sys_Filippov_representation}} is then {obtained~as
\begin{align}
& \hspace{-1mm} \dot{\bar{V}} =  \hspace{-1mm} \left\{\sum_{j \in N} \hspace{-1mm} \omega_j (-p^L_j + p_j^M - v_j - d^u_j  - \sum_{k:j\rightarrow k} p_{jk}+ \sum_{i:i\rightarrow j} p_{ij})
\right . \nonumber\\[-3mm]
& \left. \hspace{-1mm} + \hspace{-1mm} \sum_{(i,j) \in E}\hspace{-1mm} B_{ij} (\sin \eta_{ij} - \sin \eta^*_{ij}) (\omega_i - \omega_j) : \;\; v_j\in F[d^c_j(\omega_j)] \right\}.
\label{V_bar_dot}
\end{align}
}}
 {Furthermore,} from Assumption \ref{assum_passivity}, it follows that within {neighborhoods} $\Omega_j$ and $X^s_j$ of $-\omega^*_j$ and $x^{s,*}_j$ respectively it holds that
\begin{equation}\label{VSdiff}
\dot{V}^S_j(x^s_j) \leq (-\omega_j)(s_j - s^*_j) - \phi_j(-\omega_j), j \in N,
\end{equation}
noting that $\omega^* = \vect{0}_{|N|}$ follows from Lemma \ref{eqbr_lemma}.

We then consider the {function
\begin{align}\label{V_Lyapuinov}
& V(\eta, \omega, x^s) = \bar{V}(\omega,\eta) + \sum_{j \in N} V^S_j(x^s_j).
\end{align}
which as follows from \eqref{V_Lyapunov_2} and Assumption \ref{assum_passivity} is continuously differentiable with respect to the states and {has a local strict minimum at $(\eta^*, \omega^*, x^{s,*})$, as follows from Assumptions \ref{assum1} and \ref{assum_passivity},} and hence suitable to be used as a Lyapunov function {candidate}.} {Furthermore, it trivially follows that $V$ is regular, following the definition in  \cite[p. 363 - 364]{bacciotti1999stability}.}

{
{{We now consider the set-valued derivative of $V$ with respect to \eqref{sys_Filippov_representation}, i.e. $\dot{V}(x)=\left\{\frac{\partial V}{\partial x} z : z\in Q(x)\right\}$}.}
{By substituting \eqref{sys1d} and using \eqref{V_bar_dot}--\eqref{VSdiff}, it follows that}
\begin{equation}\label{e:dotV-set}
\dot{V}= \{\sum_{j \in N} (-\omega_jv_j {-\phi_j(-\omega_j)}):\;\; v_j\in F[d^c_j(\omega_j)]\}.
\end{equation}}
Using~\eqref{sys2_dc_set_valued}, we conclude that
${(\max\, \dot{V} )}\le -\sum_{j \in N} {\phi_j(-\omega_j)}\le 0,$
{where the maximum is taken over all the points in the set given by the right hand side of \eqref{e:dotV-set}.}

{By Assumptions \ref{assum1} and  \ref{assum_passivity}, the function $V$ has a strict local minimum at $x^*=(\eta^*, \omega^*, x^{s,*})$.
Then, noting that each point in $\dot V$ is nonpositive, we can construct a {compact} set $\Xi$ around the {equilibrium point, of the form $\Xi=\{(\eta, \omega, x^s) \colon V \le \epsilon, x^\ast\in \Xi, \Xi \text{ connected}\}$ {for some $\epsilon>0$} sufficiently small such that $(\eta, \omega, x^s)\in\Xi$ implies that $\omega_j, x^s_j$ lie in $\Omega_j, X^s_j$ respectively {$\forall j \in N$, {and
$\eta_{ij} < \pi/2, \forall (i,j) \in E$.}}
{Note that solutions initiated within $\Xi$ remain in $\Xi$ for all future times.
{{Therefore}, since $\Xi$ can be selected to be arbitrarily small, it follows that $x^\ast$ is Lyapunov stable.}}

Theorem 3 in \cite{bacciotti1999stability} {is now} invoked for the function $V$ on the {compact}
set $\Xi$ along solutions of \eqref{sys_Filippov_representation}.
{First, define the set $Z  =  \{(\eta, \omega, x^s) \colon 0 \in {\dot{V}\}}$ and let $M$ be the  largest weakly\footnote{{We use the notion of weakly invariant set from \cite[Dfn. 4]{bacciotti1999stability}}} invariant set within $\Xi \, \cap \, \closure{Z}$.
Then  Theorem 3 in \cite{bacciotti1999stability} guarantees that all solutions of~\eqref{sys1}--\eqref{dynsys}, \eqref{sys2_dc} with initial conditions $(\eta(0), \omega(0), x^s(0)) \in \Xi$ converge to $M$ as $t \rightarrow \infty$. }
{Note that $0 \in \dot{V}$ only if
$\omega=\omega^* = \vect{0}_{|N|}$, which also implies that {$d^{c} = \vect{0}_{|N|}$}.}
This suggests from \eqref{sys1a} and\footnote{
%{Note that $d^{c} = \vect{0}_{|N|}$ implies that the right-hand side of \eqref{sys1}--\eqref{dynsys} is Lipschitz, hence a weakly invariant set of \eqref{sys_Filippov_representation}  is also an invariant set of \eqref{sys1}--\eqref{dynsys}}.
\ak{Note that the  the right-hand side of \eqref{sys1}--\eqref{dynsys} is Lipschitz within a neighbourhood of $M$ since $\omega = \vect{0}_{|N|}$ within $M$. Hence, $M$  is also an invariant set of \eqref{sys1}--\eqref{dynsys}}.
}}
{Assumption \ref{assum_0}} {that the vectors
$\eta$ and $x^s$ are equal to some constant vectors $\bar{\eta}$ and $\bar{x}^{s}$, on the invariant set.}
Therefore, we conclude by  \cite[Thm. 3]{bacciotti1999stability} that all Filippov solutions of~\eqref{sys1}--\eqref{dynsys}, \eqref{sys2_dc} with initial conditions $(\eta(0), \omega(0), x^s(0)) \in \Xi$ converge to {the set of equilibrium points within $\Xi$.}.
\hfill $\blacksquare$

{
\begin{remark}\label{Rem:semistability}
{Theorem \ref{convthm} shows that all solutions starting  in a neighbourhood $\Xi$ of the equilibrium point $x^\ast$ will converge to the set of equilibrium points within $\Xi$. Furthermore, $\Xi$ can be chosen arbitrarily small so $x^\ast$ is also Lyapunov stable.
An extension could be obtained if {in addition to the requirements of Theorem \ref{convthm}, Assumption \ref{assum_passivity} holds} for all equilibria within $\Xi$. This means that such equilibria are all Lyapunov stable and hence arguments analogous to those in
\cite[Prop. 4.7, Thm. 4.20]{haddad2011nonlinear} can be used to deduce\footnote{{In particular, the set of equilibrium points in $\Xi$ is the $\omega$-limit set of the trajectories in $\Xi$. Furthermore, each trajectory $x(t)$ in $\Xi$ has an $\omega$-limit point  since it is bounded, i.e. it has a subsequence $x(t_n)$ that converges to an equilibrium point as $t_n\rightarrow\infty$. Hence, if the equilibrium points in $\Xi$ are Lyapunov stable one can easily deduce that $x(t)$ also converges to an equilibrium point.}} that all solutions starting in $\Xi$ converge to an equilibrium point in $\Xi$. }
\end{remark}
}

{\emph{Proof of Lemma~\ref{lemma_existence_uniq_hysteresis}:}
To show the existence of solutions, first note that for any initial condition it holds that either ${z}(0,0) \in C$ or ${z}(0,0) \in D$. {The latter results in} ${z}(0,1) \in C$ as $g(z): D \rightarrow C$. Then, from {the Lipschitz property of} the dynamics in \eqref{e:f(z)}, it follows by \cite[Proposition 2.10]{goebel2012hybrid}
that a solution to \eqref{sys4} exists. Furthermore, noting that $g: D \rightarrow C$, it follows also by \cite[Proposition 2.10]{goebel2012hybrid} that every solution to \eqref{sys4} is either complete or {it has a finite number of intervals with the last interval of the form}
$[t_{\ell},t_{{\ell}+1})\times \{{\ell}\}$ with $t_{{\ell}+1} < \infty$}}. The latter can be excluded due to the {global Lipschitz} property of the vector field in \eqref{e:f(z)}. Hence, a complete solution to \eqref{sys4} exists from any initial condition $z(0,0) \in \Lambda$.
\hfill $\blacksquare$

\emph{Proof of Proposition~\ref{dwell_time_lemma}:}
Consider any bounded solution of the system~\eqref{sys4} with states $z = (\eta, \omega, x^s, \sigma)$ and define $\epsilon_j = \omega^1_j - \omega^0_j$ following the description in \eqref{sys_hysteresis}.
{For any finite time interval between two consecutive switches {at bus $j$, i.e. $[t_{\ell,j}, t_{\ell+1,j}]$,} the value of $\dot{\omega}_j$ is bounded from above by a constant, say $d\omega^{\max}_j$. The fact that $d\omega^{\max}_j$ is finite follows from boundedness of the solution {and the fact that the vector field in \eqref{e:f(z)} is globally {Lipschitz.}}
Then,} {it follows that $t_{\ell+1,j}-t_{\ell,j}\geq \epsilon_j / d\omega^{\max}_j$.}
{Then a uniform bound between consecutive switches at any bus can be provided by $\tau = \min_{j} \epsilon_j / d\omega^{\max}_j$.}
Notice that the {bound provided in the proposition} is stated {to hold} from the second switching time to include the case  {$z(0,0) \in D$}.
\hfill $\blacksquare$

\emph{Proof of Theorem~\ref{conv_thm_hysteresis}:}
For the proof we shall make use of the {continuous} function $V$, described by \eqref{V_Lyapuinov}. Using similar arguments as in the proof of Theorem \ref{convthm} and defining $T_c =  \{t: (t,{\ell}) \in K, z(t,{\ell}) \in  {C} \}, T_d =  \{t: (t,{\ell}) \in K, z(t,{\ell}) \in  {D} \}$, {where $K$ is a hybrid time domain for \eqref{sys4} and $C$ and $D$ are defined by \eqref{e:C} and \eqref{set_D} respectively}, it follows that
\vspace{-1mm}
\begin{subequations}\label{dotV_hysteresis}
\begin{gather}
\dot{V} \leq -\sum_{j \in N} {\phi_j(-\omega_j)},  t \in T_c \\
V(g(z)) - V(z) = 0, t \in T_d,
\end{gather}
\end{subequations}
along any solution of \eqref{sys4}.  Note that when $z \in {D}$, the value of $V$ remains constant as it only depends on $x$ that is constant from \eqref{sys4_g}.

{
{Note that $V$ is a function of $x$ only, and has a strict minimum at the equilibrium point, {since $\bar{V}(\eta, \omega)$  has a minimum at $(\eta^*, 0)$ from  \eqref{V_Lyapunov_2} and $V^S_j(x^s_j)$ have strict local minima at $x^{s,*}_j$ for each $j \in N$ from Assumption \ref{assum_passivity}}. Hence} there exists a {compact}
set
$
S=\{(x, \sigma): x\in \Xi {\text{ and }} \sigma\in {\mathcal{I}}(\omega)\}$ for some {neighborhood} $\Xi$ of $x^\ast$, {such that solutions initiated in $S$ stay in $S$ for all future times}.
{The
{set} $\Xi$ is obtained in the same vein as in the proof {of~Theorem~\ref{convthm}, i.e., {$\Xi=\{(\eta, \omega, x^s) \colon V \le \epsilon, x^\ast\in \Xi, \Xi \text{ connected}\}$ {for some $\epsilon>0$} sufficiently small such that $(\eta, \omega, x^s)\in\Xi$ implies that $\omega_j, x^s_j$ lie in $\Omega_j, X^s_j$ respectively {$\forall j \in N$, {and
$\eta_{ij} < \pi/2, \forall(i,j) \in E$.}}}
Note that $\Xi$ is compact and hence from that and \eqref{dotV_hysteresis} it follows that all solutions of \eqref{sys4} that start within $\Xi$ are bounded. {Furthermore, from Lemma \ref{lemma_existence_uniq_hysteresis}  it follows that all solutions within $\Xi$ are complete.}
{In order to prove Theorem~\ref{conv_thm_hysteresis}  we make use of \cite[Corollary 8.7 (b)]{goebel2012hybrid}, which {follows from} \cite[Theorem 8.2]{goebel2012hybrid}. In particular, this Corollary can be used since from Proposition \ref{dwell_time_lemma} {we have} that the time interval between any two consecutive switches is bounded from below by a positive number. From \cite[Corollary 8.7 (b)]{goebel2012hybrid} we deduce that all complete and bounded solutions within $S$ converge to the largest weakly invariant  subset\footnote{{We use the notion of {weakly} invariant sets {provided in \cite[Dfn 6.19]{goebel2012hybrid}.}}} of the set $\{ z : V(z) = r\} \cap  S \cap \bar{u_c^{-1}(0)}$, for some $r>0$, where $u_c^{-1}(0) = \{z \in C : \dot{V} = 0\}$.}
{
{We now have that if
$\dot{V} = 0$, then} $\omega = 0$ which implies that $\sigma = 0$ by \eqref{sys_hysteresis}.
{Hence,}
{using also Assumptions \ref{assum_0} and  \eqref{sys1a}} we deduce that} $z = (\eta, \omega, x^s, \sigma)$ converges to {the set of equilibrium points in {$S$. Note that the set of equilibria of \eqref{sys4} within $S$ have continuous state $x$ that lies in  $\Xi$.}}
\hfill $\blacksquare$

{
\begin{remark}\label{rem_well_posed}
Well posedness of \eqref{sys4},  see \cite[Theorem 6.8]{goebel2012hybrid}, is required to use the invariance principle from  \cite[Thm. 8.2]{goebel2012hybrid}.
The fact that \eqref{sys4} is well posed follows trivially since the sets $C$ and $D$ are closed, the function $f$ is Lipschitz continuous
 and $g$ is \ak{outer semicontinuous}\footnote{\ak{The definition of an outer semicontinuous function is provided in \cite[Definition 5.9]{goebel2012hybrid}.}} and locally bounded \icl{relative to $D$}.
{Note also that the sets $C$ and $D$ are closed despite the discontinuity of function ${\rm sgn}(\omega_j)$ at {$\omega_j=0$} since in the definitions of $\mathcal{I}_j(\omega_j)$ and $\mathcal{I}^D_j(\omega_j)${, function} ${\rm sgn}(\omega_j)$ is used for values of $\omega_j$ that do not include $\omega=0$.}
\end{remark}
}

\bibliography{andreas_bib}

\end{document}